\documentclass[reqno,12pt]{amsart}
\usepackage{amsfonts}
\usepackage{amssymb}
\usepackage{bbm}
\usepackage{times}

\def\B{\mathcal{B}}
\def\D{\mathbb{D}}

\def\C{\mathbb{C}}

\def\N{\mathbb N}

\def\f{\frac}
\def\vect{\boldsymbol}

\def\B{\mathbb{B}}

\def\a{\alpha}

\def\msk{\medskip}

\def\bege{\begin{equation}} \def\ende{\end{equation}}
  \def\a{\alpha} 
\def\b{\beta}   
 
\def\begr{\begin{eqnarray}} \def\endr{\end{eqnarray}}

\def\bege{\begin{equation}} \def\ende{\end{equation}}
\def\begr{\begin{eqnarray}} \def\endr{\end{eqnarray}}
\def\bnum{\begin{enumerate}} \def\enum{\end{enumerate}}
\begin{document}

\title[Difference operators between weighted Bergman spaces]{Difference of composition operators between weighted Bergman spaces on the unit ball}
\author{Yecheng Shi,  Songxiao Li$^*$ and Juntao Du}

\address{Yecheng Shi\\    School of Mathematics and Statistics, Lingnan Normal University,
     Zhanjiang 524048, Guangdong, P. R. China}\email{ 09ycshi@sina.cn}

\address{Songxiao Li\\ Institute of Fundamental and Frontier Sciences, University of Electronic Science and Technology of China,
610054, Chengdu, Sichuan, P. R. China.   } \email{jyulsx@163.com}

\address{Juntao Du \\
Institute of Systems Engineering, Macau University of Science and Technology, Avenida Wai Long, Taipa, Macau. } \email{jtdu007@163.com}

\subjclass[2000]{32A36, 47B33 }
\begin{abstract}   We obtain some estimates for norm and essential norm of the difference of two composition operators between weighted Bergman spaces $A^p_\a$ and $A^q_\b$ on the unit ball. In particular, we completely characterize the boundedness and compactness of
$C_\varphi-C_\psi: A^p_\alpha\to A^q_\beta$ for full range $0<p,q<\infty$, $-1<\alpha,\beta<\infty$.
\thanks{*Corresponding author.}
\vskip 3mm \noindent{\it Keywords}: Bergman space, composition operator, difference, norm, essential norm.
\thanks{This  project was partially supported by the Macao Science and Technology Development Fund (No.186/2017/A3) and NNSF of China (No. 11720101003).   }
\end{abstract}

\maketitle

\section{Introduction}
Let $\mathbb{B}=\mathbb{B}_n$ be the open unit ball in $\C^n$. For any two points $z=(z_1,...,z_n)$ and $w=(w_1,...,w_n)$ in $\C^n$,  we write
$$\langle z,w\rangle=z_1\overline{w_1}+...+z_n\overline{w_n},$$
and
$$|z|=\sqrt{\langle z,z\rangle}=\sqrt{|z_1|^2+...+|z_n|^2}.$$
Let $H(\mathbb{\B})$ be the class of holomorphic functions on $\B$.
 Let $\varphi$ be a holomorphic self-map of $\B$. The map $\varphi$ induces a composition operator $C_\varphi$ on $H(\B)$, which is defined by
$C_\varphi f=f\circ\varphi$. We refer to \cite{cm,Sh1} for various aspects on the theory of composition operators acting on holomorphic function spaces.

Let $\rho(z,w)$ be the pseudo-hyperbolic distance between $z,w\in\B$.
Given two holomorphic self-maps $\varphi, \psi$ of $\B$, we put
$$\rho(z)=\rho(\varphi(z),\psi(z))$$
for short.  Given $\a,\beta >-1$ and $0<p,q<\infty$,   the joint pull-back measure
$ \omega_{\b,q,\varphi,\psi}$ is defined by (see \cite{KW})
$$\omega_{\b,q,\varphi,\psi}(E)=\int_{\varphi^{-1}(E)}\rho(z)^qd\nu_\b(z)+\int_{\varphi^{-1}(E)}\rho(z)^qd\nu_\b(z)$$
for Borel sets $E\subset\B$. For the simplicity, we denote $ \omega_{\b,q, \varphi,\psi}$ by $\omega_{\b,q}$.   So, $\omega_{\b,q}$ is actually the sum of two pull-back measures
$\rho^qd\nu_\b\circ\varphi^{-1}$ and $\rho^qd\nu_\b\circ\psi^{-1}$. By a standard argument one can verify that
\begr
\int_{\B}gd\omega_{\b,q}=\int_{\B}(g\circ\varphi+g\circ\psi)\rho^q d\nu_{\b}\label{gs}
\endr
for any positive Borel function $g$ on $\B$.

Let $d\nu$ be the normalized volume measure on $\B$.  For $\a>-1$,  put
$$d\nu_\a=c_\a(1-|z|^2)^\a d\nu(z),  $$
where the constant $c_\a=\f{\Gamma(n+1+\a)}{n!\Gamma(\a+1)}$ is chosen so that $\nu_\a(\B)=1$. For $0<p<\infty$ and $\a>-1$, the weighted Bergman space $A^p_\a =A^p_\a(\B)$ is the space of all $f\in H(\B)$ such that
$$\|f\|_{A^p_\a}^p=\int_{\B}|f(z)|^pd\nu_\a(z)<\infty.$$
We will also let $L^p_\a= L^p_\a(\B)$ denote the standard Lebesgue space on $\B$ with respect to the measure $\nu_\a$. The space
$A^p_\a$ equipped with the norm $\|\cdot\|_{A^p_\a}$ is a Banach space for $1\leq p<\infty$ and a complete metric space for $0<p<1$ with respect to the translation-invariant metric $(f,g)\mapsto\|f-g\|_{A^p_\a}^p$.

It seems better to clarify the concept of compact operators, since the weighted Bergman space $A^p_\a$ are not Banach space when $0<p<1$. Suppose $X$ and $Y$ are topologies vector spaces whose topologies induced by complete metrics. For a linear operator $T:X\to Y$, we   denote
$$\|T\|_{X\to Y}=\sup_{\|f\|_{X}\leq1}\|Tf\|_{Y}$$
the operator norm of $T$, where $\|\cdot\|_{X}$ and $\|\cdot\|_{Y}$ denote the norm or quasi-norm of $X$ and $Y$. If $\|T\|_{X\to Y}$ is finite, we say that $T$ is a bounded operator from $X$ to $Y$. A linear operator $T:X\to Y$ is said to be compact  if the image of every bounded sequence in $X$ has a subsequence that converges in $Y$. If $T:X\to Y$ is a bounded linear operator, then the essential norm of the operator $T:X\to Y$, denote by $\|T\|_{e,X\to Y}$, is defined as
 $$\|T\|_{e,X\to Y}=\inf\{\|T-K\|_{X\to Y}:K~ \mbox{~is compact from~} X \mbox{~to~} Y\}.$$ It is obvious that the  operator $T$ is compact if and only if $\|T\|_{e,X\to Y}=0$.

Efforts to understand the topological structure of the space of composition operators in the operator norm topology have led to the study of the operator $C_\varphi-C_\psi$ of two composition operators induced by holomorphic self-maps $\varphi, \psi$ of $\B$. In the setting of the unit disk $\D=\B_1$, by Littlewood's Subordination Principle, all composition operators, and hence all differences of two composition operators, are bounded on the Hardy space $H^p(\D)$ and the weighted Bergman space $A^p_\a(\D)$.
Thus the question of when the operator $C_\varphi-C_\psi$  is compact naturally arises.
Shapiro and Sundberg \cite{SS}   raised and studied such a question on the Hardy space, motivated by the isolation phenomenon observed by Berkson \cite{Be}.
 After that, such related problems have been studied between several spaces of analytic
functions by many authors.  See, for example,  \cite{G,NS,SL} on Hardy spaces and \cite{CKP,KW,LS, Mo,S1,S2} on weighted Bergman spaces.

In 2005, Moorhouse \cite{Mo} characterized the compact difference of composition operators on the standard weighted Bergman space $A^2_\a(\D)$ by angular derivative cancellation property. For $0<p\leq q<\infty$, Saukko  \cite{S1} obtained some compactness criterion for difference $C_\varphi-C_\psi$ from $A^p_\a(\D)$ to $A^q_\b(\D)$. In \cite{KW}, Koo and Wang gave some characterizations for the boundedness and compactness of the difference of composition operators $C_\varphi-C_\psi:A^p_\a\to A^p_\a$ on the unit ball.
  It is worth pointing out that the approach in \cite[Theorem4.5(i)]{S1} does not work as well for the unit ball. For the essential norm estimate of  $C_\varphi-C_\psi:A^p_\a\to A^q_\b$, Saukko's approach \cite{S1} is only valid for $1<p\leq q<\infty$. The main idea consists  in approximating the identity operator by a sequence of compact operators, which fits well with the study of the difference of composition operators on these spaces. In \cite{S1}, this sequence was the finite rank operators $S_n$ which map $A^p_\a$ to $n$th-partial sum of the Talor series of $f$. However, this sequence is uniformly bounded only for $p>1$. The first goal of this paper is to study the norm and the essential norm of the difference of composition operators  $C_\varphi-C_\psi:A^p_\a\to A^q_\b$  for   $0<p\leq q<\infty$. Our work requires certain new approach and substantial amount of extra works.

The first two main results of this paper are the following theorems. For the simplicity, we denote
$$
\Gamma(\varphi, \psi )=\Big(\frac{(1-|a|^2)^{s}}{|1-\langle a,\varphi(z)\rangle|^{(n+1+\alpha)\lambda+s}}+\frac{(1-|a|^2)^{s}}{|1-\langle a,\psi(z)\rangle|^{(n+1+\alpha)\lambda+s}}\Big) .
$$ 

\noindent{\bf Theorem 1.1.}  {\it Let $0<p\leq q< \infty$, $-1<\alpha, \beta<\infty$. Suppose $\varphi$ and $\psi$ are holomorphic self-maps of $\B$, and denote $\lambda=q/p$. Then the  operator $C_\varphi-C_\psi$ maps $A^p_\alpha$ into $A^q_\beta$
if and only if the joint pull-back
measure $\omega_{\b,q}$ is a $(\lambda,\a)$-Carleson measure. Furthermore,
\begr
 \|C_\varphi-C_\psi\|_{A^p_\a\to A^q_\b}^q \asymp \sup_{a\in\B}\int_{\B}    \Gamma(\varphi, \psi )     \rho(z)^qd\nu_\beta(z)\nonumber
\endr
for some (or equivalent for all)  $s>0$.
}\msk

\noindent{\bf Theorem 1.2.}  {\it Let $0<p\leq q< \infty$, $-1<\alpha, \beta<\infty$. Suppose $\varphi$ and $\psi$ are holomorphic self-maps of $\B$. Then the operator $C_\varphi-C_\psi$ maps $A^p_\alpha$ into $A^q_\beta$
if and only if the joint pull-back measure $\omega_{\b,q}$ is a vanishing $(\lambda,\a)$-Carleson measure.
Furthermore,
\begr
 \|C_\varphi-C_\psi\|_{e,A^p_\a\to A^q_\b}^q \asymp \limsup_{|a|\to1}\int_{\B}\Gamma(\varphi, \psi )\rho(z)^qd\nu_\beta(z),  \nonumber
\endr
for some (or equivalent for all)  $s>0$.  }\msk

  Saukko in \cite{S2} characterized the bounded difference of composition operators from $A^p_\a(\D)$ into Lebesgue spaces $L^q(\D,\mu)$ when $\a>-1$ and $p>q$. In particular, the following result was shown in \cite{S2}.\msk

\noindent{\bf Theorem A.} {\it Let $0<q<p<\infty$, $\a>-1$ and $\mu$ a positive Borel measure on $\D$. Denote $s:=p/(p-q)$. Let $\varphi$ and $\psi$ be analytic self-maps of $\D$, and denote $\sigma(z)=\left|\frac{\varphi(z)-\psi(z)}{1-\overline{\varphi(z)}\psi(z)}\right|$ for every $z\in\D$. Then the following are equivalent:

(i) the operator $C_\varphi-C_\psi$ maps $A^p_\a(\D)$ into $L^q(\D,\mu)$;

(ii) the operators $\sigma C_\varphi$ and $\sigma C_\psi$ map $A^p_\a(\D)$ into $L^q(\D,\mu)$;

(iii) the function
$$K_{\varphi,\psi}(z):=\int_{\D}\left|\left(\frac{1-|z|^2}{1-\overline{z}\varphi(w)}\right)^{\frac{\a+2}{p}}-\left(\frac{1-|z|^2}{1-\overline{z}\varphi(w)}\right)^{\frac{\a+2}{p}}\right|^qd\mu(w)$$
belongs to $L^s(\D,A_\a)$.  }\msk

The last main result of this paper (Theorem 1.3)   gives an estimate for the norm of the difference of composition operators $C_\varphi-C_\psi:A^p_\a\to A^q_\b$, when $0<q<p<\infty$. Saukko \cite{S2}   pointed out that, by Pitt's theorem, the  operator $C_\varphi-C_\psi:A^p_\a\to L^q(\mu)$ is compact, whenever it is bounded, for $1\leq q<p<\infty$. In fact, we can prove that this phenomenon is true for full range $0<q<p<\infty$. Our method is new even for the   case of the unit disk. \msk

\noindent{\bf Theorem 1.3.}  {\it Let $0<q<p< \infty$, $-1<\alpha, \beta<\infty$. Suppose $\varphi$ and $\psi$ are holomorphic self-maps of $\B$. Set $t=\frac{p}{p-q}$.
Then the following are equivalent:

(i)the operator  $C_\varphi-C_\psi: A^p_\alpha\to A^q_\beta$ is bounded;

(ii) the operator $C_\varphi-C_\psi: A^p_\alpha\to A^q_\beta$ is compact;

(iii) the joint pull-back measure  $\omega_{\b,q}$ is a $(\lambda,\a)$-Carleson measure.
Furthermore,
\begr
&&\|C_\varphi-C_\psi\|_{A^p_\a\to A^q_\b}^q\nonumber\\
&\asymp&\left\|\int_{\B}\left(\frac{(1-|a|^2)^{s}}{|1-\langle a,\varphi(z)\rangle|^{n+1+\alpha+s}}+\frac{(1-|a|^2)^{s}}{|1-\langle a,\psi(z)\rangle|^{n+1+\alpha+s}}\right)\rho(z)^qd\nu_\beta(z)\right\|_{L^{t}(\B,\nu_\a)}, \nonumber
\endr
for some (or equivalent for all)  $s>0$.  }\msk

The present paper is organized as follows. In Section 2, we give some notations and preliminary results which will be used later. In Sections 3, we give the proofs for Theorems 1.1, 1.2 and 1.3.  For two quantities $A$ and $B$, we use the abbreviation
$A\lesssim B$ whenever there is a positive constant $C$ (independent of the associated variables) such that $A\leq CB$.
We write $A\asymp B$, if $A\lesssim B\lesssim A$.

\section{prerequisites}\vspace{0.2truecm}

In this section we introduce some notations and recall some well known results that will be used throughout the paper.

\subsection{Pseudo-hyperbolic distance}

For any $z\in\B$, let $P_z$ be the orthogonal projection from $\C^n$ onto the one dimensional subspace $[z]=\{\lambda z:\lambda\in\C\}$ generated by $z$, and $Q_z$ be the orthogonal projection from $\C^n$ onto $\C^n\ominus[z]$.
Thus $P_0(w)=0$,  $Q_0(w)=w$ and
$$P_z(w)=\f{\langle w,z\rangle}{|z|^2}z, ~~~~\,~~~\,~~Q_z(w)=w-\f{\langle w,z\rangle}{|z|^2}z,\mbox{~~~if~~}z\neq0.$$
We denote by $\sigma_z(w)$ the M\"{o}bius transformation on $\B$ that interchanges the points $0$ and $z$. More explicitly,
$$\sigma_z(w)=\f{z-P_z(w)-\sqrt{1-|z|^2}Q_z(w)}{1-\langle w,z\rangle}.$$
Note that $P_z(w)=w$ when $n=1$.   It is well known that $\sigma_z$ satisfies the following properties:
$\sigma_z\circ\sigma_z(w)=w$, and
$$1-|\sigma_z(w)|^2=\frac{(1-|z|^2)(1-|w|^2)}{|1-\langle w,z\rangle|^2}, ~~~~z,w\in\B.$$
For $z,w\in\B$, the pseudo-hyperbolic distance between $z$ and $w$ is defined by
$$\rho(z,w)=|\sigma_z(w)|,$$
while the Bergman metric   is given by
 $$\b(z,w)=\frac{1}{2}\log\frac{1+\rho(z,w)}{1-\rho(z,w)}.$$
 It is also well known that the pseudo-hyperbolic metric have the following strong form of triangle inequality (see \cite{DW}):
 $$\rho(z,w)\leq\f{\rho(z,a)+\rho(a,w)}{1+\rho(z,a)\rho(a,w)}$$
for all $a,z,w\in\B$.  For $z\in \B$ and $r>0$, the Bergman metric  ball at $z$ is denoted by
 $$D(z,r)=\{w\in\B:\b(z,w)<r\},$$
and the pseudo-hyperbolic ball at $z\in\B$ with radius $r\in(0,1)$ is given by
 $$\triangle(z,r)=\{w\in\B:\rho(z,w)<r\}.$$
 Note that $\rho(z,0)=|z|$ since $\sigma_0(z)=-z$, so $\triangle(0,r)$ is a Euclidean ball $|z|<r$.  For any $z\neq0$, $\triangle(z,r)$ is an ellipsoid consisting of all $w\in \B$ such that
$$\f{|P_z(w)-c|^2}{r^2t^2}+\f{|Q_z(w)|^2}{r^2t}<1,$$
where
$$ c=\f{(1-r^2)z}{1-r^2|z|^2}~~,\mbox{~~~~}~~t=\f{1-|z|^2}{1-r^2|z|^2}.$$
Furthermore, if $0<r<1$, then the weighted volume
 $$\nu_\a(\triangle(z,r))\asymp (1-|z|^2)^{n+1+\a}.$$

The following lemma should be known to some experts, but we cannot find a reference.
So we give the proof for completeness.\msk

\noindent{\bf Lemma 2.1.}  {\it The pseudo-hyperbolic metric
\begr
\rho(z,w)\leq\left|\f{z-w}{1-\langle z,w\rangle}\right|,\nonumber
\endr
for all $z,w\in\B$.
}\msk

\noindent{\bf Proof.}
Since
\begr
&&|z-P_z(w)-\sqrt{1-|z|^2}Q_z(w)|^2\nonumber\\
&=&|P_z(z-w)+\sqrt{1-|z|^2}Q_z(z-w)|^2\nonumber\\
&=&|P_z(z-w)|^2+(1-|z|^2)|Q_z(z-w)|^2\nonumber\\
&=&(1-|z|^2)|z-w|^2+|\langle z-w,z\rangle|^2 \leq |z-w|^2,\nonumber
\endr
we get the desired result. $\Box$

\subsection{Local estimates and test functions}\vspace{0.3truecm}

The following lemmas are crucial in our work and will be repeatedly used throughout the paper.\msk

\noindent{\bf Lemma 2.2.}  {\it Let $0<p\leq q<\infty$, $\alpha>-1$ and $0<s<r<1$ be arbitrary. Then there exists a constant $C=C(p,q,\alpha,s,r)$ such that
\begr
|f(z)-f(a)|^q\leq C\rho(z,a)^q\frac{\int_{\bigtriangleup(a,r)}|f(w)|^pd\nu_\alpha(w)}{(1-|a|^2)^{(n+1+\alpha)q/p}}\nonumber
\endr
for all $a\in\B$, $z\in \bigtriangleup(a,s)$ and $f\in A^p_\alpha$ with $\|f\|_{A^p_\alpha}\leq1$.
}\msk

\noindent{\bf Proof.} For the case $p=q$ see \cite[Lemma 2.2]{KW}. The case $p>q$ can be proved similarly as \cite[Lemma 3.1]{S1}. $\Box$\msk

Let $\vect{e_j}=(0,\cdots,0,1,0,\cdots,0)\in \C^n$, where $1\leq j\leq n$ and $1$ is on the $j$-th component. For
$0<t<1$, let
$$\vect{t_1}=t\vect{e_1},\mbox{~~~~}\vect{t_j}=t^2\vect{e_1}+t\sqrt{1-t^2}\vect{e_j}~~\mbox{~~~}(j=2,..,n).$$
For $N>0$ and $0<t<1$, let
$$t_N=1-N(1-t).$$

We need the following two results from \cite{KW}.\msk

\noindent{\bf Lemma 2.3.}  {\it Suppose $s>1$ and $0<r_0<1$. Then there are $N=N(r_0)>1$ and $C=C(s,r_0)$ such that
\begr
&&\left|\frac{1}{(1-\langle a,\vect{t_1}\rangle)^s}-\f{1}{(1-\langle b,\vect{t_1}\rangle)^s}\right|+
\sum_{j=1}^n\left|\frac{1}{(1-\langle a,t_N\vect{t_j}\rangle)^s}-\f{1}{(1-\langle b,t_N\vect{t_j}\rangle)^s}\right|\nonumber\\
&\geq& C\rho(a,b)\left|\frac{1}{(1-\langle a,\vect{t_1}\rangle)^s}\right|, \nonumber
\endr
for all $a\in\triangle(t\vect{e_1},r_0)$ with $1-t<\f{1}{2N}$ and $b\in\B$.
}\msk

\noindent{\bf Lemma 2.4.}  {\it Suppose $0<p<\infty$, $\a>-1$, $a\in\B$ and $0<r_0<1$.
Let $\delta>0$ such that $t=n+1+\a+\delta>p$. Let $N=N(r_0)$ be as in Lemma 2.3, $|a|_N=1-N(1-|a|)$
and
$$\vect{|a|_1}=|a|\vect{e_1},\mbox{~~~~}\vect{|a|_j}=|a|^2\vect{e_1}+|a|\sqrt{1-|a|^2}\vect{e_j}~~\mbox{~~~}(j=2,..,n).$$
Let
$$f_{a,0}(z)=\frac{(1-|a|^2)^{\delta/p}}{(1-\langle z,\sigma^{-1}(\vect{|a|_1}))^{t/p}}$$
and
$$f_{a,j}(z)=\frac{(1-|a|^2)^{\delta/p}}{(1-|a|_N\langle z,\sigma^{-1}(\vect{|a|_j})\rangle)^{t/p}}~~\mbox{~~~~}~~(j=1,..,n),$$
where $\sigma$ is a rotation which maps $a$ to $|a|\vect{e_1}$.
Then $\sum\limits_{j=0}^n\|f_{a,j}\|_{A^p_\a}\lesssim1$, and
$$\sum_{j=0}^n|f_{a,j}(z)-f_{a,j}(w)|\geq C(\a,r_0)\rho(z,w)|f_{a,0}(z)|$$
for any $z\in\triangle(a,r_0)$, $w\in\B$ with $|a|>1-\f{1}{2N}$.}\msk

\noindent{\bf Remark.}  { The above lemma can be found in the proof of Theorem 3.1 of \cite{KW}. }

\subsection{Carleson measure}\vspace{0.3truecm} Let $\mu$ be a positive Borel measure on $\B$. For $\lambda>0$ and $\a>-1$, we say that $\mu$ is a $(\lambda,\a)$-Bergman Carleson measure if for any two positive numbers
$p$ and $q$ with $q/p=\lambda$ there is a positive constant $C>0$ such that
$$\int_{\B}|f(z)|^qd\mu(z)\leq C\|f\|_{A^p_\a}^q$$
for any $f\in A^p_\a$. We also denote by
$$\|\mu\|_{\lambda,\a}=\sup_{f\in A^p_\a,\|f\|_{A^p_\a}\leq1}\int_{\B}|f(z)|^qd\mu(z).$$

The Bergman Carleson measure was first studied by Hastings\cite{Ha}, and independently by Oleinik and Pavolv \cite{OP} and Oleinik \cite{O},
and further pursued by Luecking\cite{Lu1,Lu2}, Cima and Wogen\cite{CW}, and many others. The statement in terms of pseudo-hyperbolic balls, essentially due to Luecking \cite{Lu1}, is more convenient to use in this paper. The following result can be found in \cite[Theorem A]{PZ}.\msk

\noindent{\bf Theorem B.}  {\it For a positive Borel measure $\mu$ on $\B$, $0<p\leq q< \infty$, $-1<\a<\infty$ and $0<r<1$, the following are equivalent:

(i) There is a constant $C_1>0$ such that for any $f\in A^p_\a$
$$\int_{\B}|f(z)|^qd\mu(z)\leq C_1\|f\|_{A^p_\a}^q.$$

(ii)  $$\|\mu\|_{\lambda,\a,r}=\sup_{a\in\B}\frac{\mu(\bigtriangleup(a,r))}{(1-|a|^2)^{(n+1+\a)q/p}}<\infty.$$

(iii) There is a constant $C_2>0$ such that, for some (every) $t>0$,
 $$\sup_{a\in\B}\int_{\B}\frac{(1-|a|^2)^t}{|1-\langle z,a\rangle|^{(n+1+\a)q/p+t}}d\mu(z)\leq C_2.$$

Furthermore, the constants $C_1$, $C_2$, and $\|\mu\|_{\lambda,\a,r}$ are all comparable to $\|\mu\|_{\lambda,\a}$ with $\lambda=q/p$.}\msk

We say that $\mu$ is a vanishing $(\lambda,\a)$-Bergman Carleson measure if for any two positive numbers $p$ and $q$ satisfying $q/p=\lambda$ and any sequence $\{f_k\}$ in $A^p_\a$ with $\|f\|_{A^p_\a}\leq1$ and $f_k(z)\to0$ uniformly on any compact subset of $\B$,
$$\lim_{k\to\infty}\int_{\B}|f_k(z)|^qd\mu(z)=0.$$
It is well known that, for $\lambda\geq1$, $\mu$ is a vanishing $(\lambda,\a)$-Bergman Carleson measure if and only if
\begr
\lim_{|a|\to1}\int_{\B}\frac{(1-|a|^2)^t}{|1-\langle z,a\rangle|^{(n+1+\a)q/p+t}}d\mu(z)=0\label{1.2}
\endr
for some (any) $t>0$.  For $s\in(0,1)$,  denote $\B_s=\{z\in\B:|z|<s\}$. Let $\mu$ be a positive Borel measure on $\B$. We denote by $\mu|_{\B\backslash\B_s}$ the restriction of $\mu$ to
$\B\backslash \B_s$.\msk

\noindent{\bf Lemma 2.5.}  {\it Let $0<p\leq q< \infty$, $-1<\a<\infty$. Suppose the positive Borel measure $\mu$ on $\B$ is a $(\lambda,\a)$-Bergman Carleson measure with $\lambda=q/p$. Then
  \begr
\limsup_{|a|\to1}\frac{\mu(\bigtriangleup(a,r))}{(1-|a|^2)^{(n+1+\a)q/p}}
&=&\limsup_{s\to1}\|\mu|_{\B\backslash\B_s}\|_{\lambda,\a,r}\nonumber\\
&\asymp&\limsup_{|a|\to1}\int_{\B}\frac{(1-|a|^2)^t}{|1-\langle z,a\rangle|^{(n+1+\a)q/p+t}}d\mu(z)\nonumber
\endr
for some (any) $t>0$ and some (any) $r\in (0,1)$.}\msk

\noindent{\bf Proof.}  We first  prove that
$$\limsup_{s\to1}\|\mu|_{\B\backslash\B_s}\|_{\lambda,\a,r}=\limsup_{|a|\to1}\frac{\mu(\bigtriangleup(a,r)))}{(1-|a|^2)^{(n+1+\a)q/p}}.$$
Let $t_r(s)=\f{s-r}{1-sr}$ and $s>r$. By the ellipsoid description of $\triangle(a,r)$, we known that the point
$c=\f{1-r^2}{1-r^2|a|^2}a$ is the center of the ellipsoid $\triangle(a,r)$, $c\in[a]$, the intersection of $\triangle(a,r)$ with $[a]$ is a one-dimensional disk of radius $\f{1-|a|^2}{1-r^2|a|^2}r$ and the intersection of $\triangle(a,r)$ with $\C^n\ominus[a]$ is an $(n-1)$-dimensional Euclidean ball of radius $r\sqrt{\f{1-|a|^2}{1-r^2|a|^2}}$. Therefore, $\triangle(a,r)\subseteq\B_s$ if and only if
$$\f{1-r^2}{1-r^2|a|^2}|a|+\f{1-|a|^2}{1-r^2|a|^2}r\leq s.$$
After a calculation, we get
$\triangle(a,r)\cap(\B\backslash\B_s)\neq \emptyset$ if and only if $|a|> t_r(s)$.
It is easy to see that $t_r(s)$ is continuous and increasing on $[r,1)$, and $\lim_{s\to1}t_r(s)=1$.
Thus,
\begr
\limsup_{|a|\to1}\frac{\mu(\bigtriangleup(a,r)))}{(1-|a|^2)^{(n+1+\a)q/p}}&=&\lim_{s\to1}\sup_{|a|\geq t_r(s)}\frac{\mu(\bigtriangleup(a,r)))}{(1-|a|^2)^{(n+1+\a)q/p}}\nonumber\\
&\geq&\lim_{s\to1}\sup_{|a|\geq t_r(s)}\frac{\mu(\bigtriangleup(a,r))\cap(\B\backslash\B_s))}{(1-|a|^2)^{(n+1+\a)q/p}}\nonumber\\
&=&\lim_{s\to1}\sup_{a\in\B}\frac{\mu(\bigtriangleup(a,r))\cap(\B\backslash\B_s))}{(1-|a|^2)^{(n+1+\a)q/p}}\nonumber\\
&=&\lim_{s\to1}\|\mu|_{\B\backslash\B_s}\|_{\lambda,\a,r}.\nonumber
\endr
On the other hand, denote $A=\limsup\limits_{s\to1}\|\mu|_{\B\backslash\B_s}\|_{\lambda,\a,r}$. For any $\epsilon>0$, there exists $0<t_*<1$, such that if $t_*\leq s<1$, we have
$$\|\mu|_{\B \backslash\B_s}\|_{\lambda,\a,r}<A+\epsilon.$$
For any fixed $s$, we know that $\bigtriangleup(a,r)\subset\B\backslash\B_s$, as $|a|$ close enough to 1. Therefore, there exists a $0<l<1$, such that
\begr
\|\mu|_{\B\backslash\B_s}\|_{\lambda,\a,r}&=&\sup_{a\in\B}\frac{\mu(\bigtriangleup(a,r)\cap(\B\backslash\B_s))}{(1-|a|^2)^{(n+1+\a)q/p}}\nonumber\\
&\geq&\sup_{|a|>l}\frac{\mu(\bigtriangleup(a,r))}{(1-|a|^2)^{(n+1+\a)q/p}}.\nonumber
\endr
Hence,
\begr
A+\epsilon&\geq&
\limsup_{|a|\to1}\frac{\mu(\bigtriangleup(a,r)))}{(1-|a|^2)^{(n+1+\a)q/p}}.\nonumber
\endr
Since $\epsilon$ is arbitrary, we have
\begr
\limsup_{s\to1}\|\mu|_{\B\backslash\B_s}\|_{\lambda,\a,r}&\geq&
\limsup_{|a|\to1}\frac{\mu(\bigtriangleup(a,r)))}{(1-|a|^2)^{(n+1+\a)q/p}}.\nonumber
\endr

If $z\in\triangle(a,r)$, then $1-|a|^2\asymp |1-\langle z,a\rangle|$. Thus, we have
\begr
\limsup_{|a|\to1}\frac{\mu(\bigtriangleup(a,r))}{(1-|a|^2)^{(n+1+\a)q/p}}
&\lesssim&\limsup_{|a|\to1}\int_{\triangle(a,r)}\frac{(1-|a|^2)^t}{|1-\langle z,a\rangle|^{(n+1+\a)q/p+t}}d\mu(z)\nonumber\\
&\leq&\limsup_{|a|\to1}\int_{\B}\frac{(1-|a|^2)^t}{|1-\langle z,a\rangle|^{(n+1+\a)q/p+t}}d\mu(z).\nonumber
\endr

For any fixed $s$, by Theorem B, we have
\begr
&&\limsup_{|a|\to1}\int_{\B}\frac{(1-|a|^2)^t}{|1-\langle z,a\rangle|^{(n+1+\a)q/p+t}}d\mu(z)\nonumber\\
&=&\limsup_{|a|\to1}\left(\int_{\B_s}+\int_{\B\backslash \B_s}\right)\frac{(1-|a|^2)^t}{|1-\langle z,a\rangle|^{(n+1+\a)q/p+t}}d\mu(z)\nonumber\\
&\leq&\limsup_{|a|\to1}\frac{(1-|a|^2)^t\mu(\B_s)}{(1-s)^{(n+1+\a)q/p+t}}+\sup_{a\in\B}\int_{\B}\frac{(1-|a|^2)^t}{|1-\langle z,a\rangle|^{(n+1+\a)q/p+t}}d\mu|_{\B\backslash \B_s}(z)\nonumber\\
&\lesssim&\|\mu|_{\B\backslash\B_s}\|_{\lambda,\a,r}.\nonumber
\endr
Letting $s$ tend to $1$, we get
\begr
\limsup_{|a|\to1}\int_{\B}\frac{(1-|a|^2)^t}{|1-\langle z,a\rangle|^{(n+1+\a)q/p+t}}d\mu(z)
\lesssim\limsup_{s\to1}\|\mu|_{\B\backslash\B_s}\|_{\lambda,\a,r}. \nonumber
\endr
The proof is complete. $\Box$ \msk

\noindent{\bf Definition 2.1.}  {A sequence $\{a_k\}$ of distinct points in $\B$ is call a
 \emph{separated} sequence in the pseudo-hyperbolic metric if $\delta_0:=\inf_{i\neq j}\rho(a_i,a_j)>0$. The number $\delta_0$ is called the \emph{ separated constant} of $\{a_k\}$. We say that the sequence $\{a_k\}$ is \emph{$\delta$-separated}, if $0<\delta\leq\delta_0$.
 }\msk

 The following lemma can be found in \cite[Lemma 5]{DW} or \cite[Lemma 3]{Lu4}\msk.

 \noindent{\bf Lemma 2.6.} {\it If $\{a_k\}$ is a separated sequence in $\B$ with separation constant $\delta_0$. For $z\in\B$ and $0<r<1$, let $L$ denote the number of points in $\{a_k\}$ that lie in the pseudohyperbolic ball $\triangle(z,r)$. Then
$$L\leq\left(\f{2}{\delta_0}+1\right)^{2n}\f{1}{(1-r^2)^n}.$$   }\msk

\noindent{\bf Definition 2.2.}  { Suppose $0<r<1$. A sequence $\{a_k\}$ of distinct points in $\B$ is call an
 \emph{$r$-lattice} in the pseudo-hyperbolic metric if it is $r$-separated and $\B=\bigcup\limits_{i=1}^{\infty}\triangle(a_k,r)$.
 }\msk

 \noindent{\bf Remark.}  {(1). In \cite{Zhu} the definition of $r$-lattice is slightly different to ours but it causes no difficulties as we have only notice that $D(a_k,r)=\triangle(a_k,\tanh(r))$.\msk

 (2). By Lemma 2.6, similarly as \cite[Lemma 4]{Lu4}, it is easy to see that there exists a $r$-lattice for any $0<r<1$.}\msk

 The following result is essentially due to Luecking (\cite{Lu4}) and can be found in \cite[Theorem B]{PZ}. \msk

\noindent{\bf Theorem C.}  {\it For a positive Borel measure $\mu$ on $\B$, $0<q<p<\infty$ and $-1<\a<\infty$, the following statements are equivalent:

(i) There is a constant $C_1>0$ such that for any $f\in A^p_\a$
$$\int_{\B}|f(z)|^qd\mu(z)\leq C_1\|f\|_{A^p_\a}^q.$$

(ii) The function
$$\hat{\mu}_r(z):=\f{\triangle(z,r)}{(1-|z|^2)^{n+1+\a}}$$
is in $L^{p/(p-q)}(\nu_\a)$ for any (some) fixed $r\in(0,1)$.

(iii) For any $r$-lattice $\{a_k\}$, the sequence
$$\{\mu_k\}:=\left\{\f{\mu(\triangle(a_k,r))}{(1-|a_k|^2)^{(n+1+\a)q/p}}\right\}$$
belongs to $l^{p/(p-q)}$ for any (some) fixed $r\in(0,1)$.

(iv) For any $s>0$, the Berezin-type transform $B_{s,\a}(\mu)$ belongs to $L^{p/(p-q)}_{\nu_\a}$.

Furthermore, with $\lambda=q/p$, one has
$$\|\hat{\mu}_r\|_{L^{p/(p-q)}_{\nu_\a}}\asymp\|\{\mu_k\}\|_{l^{p/(p-q)}}\asymp\|B_{s,\a}(\mu)\|_{L^{p/(p-q)}_{\nu_\a}}\asymp\|\mu\|_{\lambda,\a}.$$
Here, for a positive measure $\mu$, the Berezin-type transform $B_{s,\a}(\mu)$ is
$$B_{s,\a}(\mu)(z)=\int_{\B}\f{(1-|z|^2)^s}{|1-\langle z,w\rangle|^{n+1+\a+s}}d\mu(w).$$
}

It is well known that, when $0<q<p<\infty$, that is $0<\lambda=q/p<1$, $\mu$ is a vanishing $(\lambda,\a)$-Bergman Carleson measure if and only if $\mu$ is a $(\lambda,\a)$-Bergman Carleson measure (see \cite{ZZ}).\msk

\section{Proof of the main results}

\noindent{\bf Proof of Theorem 1.1.} We first give the upper estimates for the norm of the   operator $C_\varphi-C_\psi$.
Fix $0<s_0<r_0<1$, set $E=\{z\in\B:\rho(z)\geq s_0\}$ and $E^\prime=\B\backslash E$. Let $f\in A^p_\a$ with $\|f\|_{A^p_\a}\leq1$. Then
$$\|(C_\varphi-C_\psi)f\|_{A^q_\b}^q=\Big(\int_{E}+\int_{E^\prime}\Big)  |f\circ\varphi(z)-f\circ\psi(z)|^qd\nu_\b(z).$$
Using (\ref{gs}), the first term is uniformly bounded above by
$$\big(\f{2}{s_0}\big)^q\int_{E}|f(z)|^qd\omega_{\b,q}(z)\lesssim \|\omega_{\b,q}\|_{\lambda,\a}^q.$$
Applying Lemma 2.2, Fubini's Theorem and $1-|z|^2\asymp1-|w|^2$ for $z\in\triangle(w,r_0)$, we see that the second term is bounded  by
\begr
&&C\int_{\B}|f(w)|^p \f{\int_{\varphi^{-1}(\triangle(w,r_0))\cap \{z\in\B:\rho(z)<s_0\}}\rho(z)^qd\nu_\b(z)}{(1-|w|^2)^{(n+1+\a)q/p}}d\nu_\a(w)\nonumber\\
&\lesssim&\int_{\B}|f(w)|^p \f{\omega_{\b,q}(\triangle(w,r_0))}{(1-|w|^2)^{(n+1+\a)q/p}}d\nu_\a(w)\nonumber\\
&\lesssim&\|\omega_{\b,q}\|_{\lambda,\a,r_0}.  \nonumber
\endr

Next, we give the lower estimate for the norm of the   operator $C_\varphi-C_\psi$. Denote $r_1=1-\f{1}{2N}$, where $N$ is defined as in Lemma 2.3. By Lemma 2.4, we have
\begr
\|C_\varphi-C_\psi\|_{A^p_\a\to A^q_\b}^q&\gtrsim&\sup_{|a|> r_1}\sum_{j=0}^n\|(C_\varphi-C_\psi)f_{a,j}\|_{A^q_\b}^q\nonumber\\
&\gtrsim& \sup_{|a|> r_1}\int_{\varphi^{-1}(\triangle(a,r_0))}\f{\rho(z)^q}{(1-|a|^2)^{(n+1+\a)q/p}}d\nu_\b(z),  \nonumber
\endr
where we used the fact $|1-\langle z,a\rangle|\asymp1-|a|$ for $z\in \triangle(a,r_0)$.  Similarly
\begr
\|C_\varphi-C_\psi\|_{A^p_\a\to A^q_\b}^q
\gtrsim \sup_{|a|> r_1}\int_{\psi^{-1}(\triangle(a,r_0))}\f{\rho(z)^q}{(1-|a|^2)^{(n+1+\a)q/p}}d\nu_\b(z).\nonumber
\endr
Therefore,
\begr
\|C_\varphi-C_\psi\|_{A^p_\a\to A^q_\b}^q
\gtrsim \sup_{|a|>r_1}\f{\omega_{\b,q}(\triangle(a,r_0))}{(1-|a|^2)^{(n+1+\a)q/p}}.\nonumber
\endr
For $|a|\leq r_1$, take $r_2=\f{r_0+r_1}{1+r_0r_1}$, then
$\triangle(a,r_0)\subset\triangle(0,r_2)$. Therefore, by Lemma 2.1,
\begr
&&\f{\omega_{\b,q}(\triangle(a,r_0))}{(1-|a|^2)^{(n+1+\a)q/p}}\nonumber\\
&\leq& \f{1}{(1-r_1^2)^{(n+1+\a)q/p}} \left(\int_{\varphi^{-1}(\triangle(a,r_0))}\rho(z)^qd\nu_\b(z)+\int_{\psi^{-1}(\triangle(a,r_0))}\rho(z)^qd\nu_\b(z)\right)\nonumber\\
&\leq& \f{1}{(1-r_1^2)^{(n+1+\a)q/p}}\left(\int_{\varphi^{-1}(\triangle(a,r_0))}+\int_{\psi^{-1}(\triangle(a,r_0))}\right)\f{|\varphi(z)-\psi(z)|^q}{(1-r_2)^q}d\nu_\b(z)\nonumber\\
&\lesssim& \|C_\varphi-C_\psi\|_{A^p_\a\to A^q_\b}^q.\nonumber%\cdot\|z\|_{A^p_\a}^q.\nonumber
\endr\msk
Then the result follows by Theorem B and (\ref{gs}). The proof is complete.   $\Box$

\noindent{\bf Lemma 3.1.}  {\it Let $0<p\leq q< \infty$, $-1<\a,\b<\infty$. There is a constant $C>0$ such that
\begr
\|C_\varphi-C_\psi\|_{e,A^p_\a\to A^q_\b}^q\geq C\limsup_{|a|\to1}\sum_{j=0}^n\|(C_\varphi-C_\psi)f_{a,j}\|_{A^q_\b}^q.\nonumber
\endr
Here $f_{a,j}$ is defined in Lemma 2.4. }\msk

\noindent{\bf Proof.}  Let $K$ be a compact operator from $A^p_\a$ into $A^q_\b$. Consider the operator on $H(\B)$ defined by
$$K_m(f)(z)=f(\f{m}{m+1}z),~~m\in\N.$$
Denote $R_m=I-K_m$. It is easy to see that $K_m$ is compact on $A^q_\b$ and 
$$\|K_m\|_{A^q_\b\to A^q_\b}\leq1, \,~~~~\,  \|R_m\|_{A^q_\b\to A^q_\b}\leq2$$ for any positive integer $m$.
Then we have
\begr
2\|(C_\varphi-C_\psi)-K\|_{A^p_\a\to A^q_\b} &\geq& \|R_m\circ(C_\varphi-C_\psi-K)\|_{A^p_\a\to A^q_\b}\nonumber\\
&\gtrsim&\sup_{a\in\B}\|R_m\circ(C_\varphi-C_\psi-K)(f_{a,j})\|_{A^q_\b}.\nonumber
\endr
Since $K$ is compact, we can extract a sequence $\{a_i\}\subset\B$ such that $|a_i|\to1$ and $Kf_{a_i,j}$ converges to some $f_j\in A^q_\b$ for $j=0,1,\cdots,n$.
So, when $0<q<1$, we get
\begr
&&\|R_m\circ(C_\varphi-C_\psi-K)(f_{a_i,j})\|_{A^q_\b}^q\nonumber\\
&\geq&\|R_m\circ(C_\varphi-C_\psi)(f_{a_i,j})\|_{A^q_\b}^q
-\|R_m\circ K(f_{a_i,j})\|_{A^q_\b}^q\nonumber\\
&\geq&\|(C_\varphi-C_\psi)(f_{a_i,j})\|_{A^q_\b}^q-\|K_m\circ(C_\varphi-C_\psi)(f_{a_i,j})\|_{A^q_\b}^q\nonumber\\
&~~&-\|R_m(K(f_{a_i,j})-f_j)\|_{A^q_\b}^q-\|R_m(f_j)\|_{A^q_\b}^q.\label{1}
\endr
Since $\varphi(\f{m}{m+1}\B)$ and $\psi(\f{m}{m+1}\B)$ are contained in a compact subset of $\B$, the Cauchy-Schwartz inequality yields, for every $z\in\B$
\begr
|K_m\circ(C_\varphi-C_\psi)(f_{a,0})(z)|
&\leq& \f{(1-|a|^2)^{\delta/p}}{|1-\langle \varphi(\f{m}{m+1}z),\sigma^{-1}(|a|\vect{e_1})\rangle|^{t/p}}\nonumber\\
&+&\f{(1-|a|^2)^{\delta/p}}{|1-\langle \psi(\f{m}{m+1}z),\sigma^{-1}(|a|\vect{e_1})\rangle|^{t/p}}\nonumber\\
&\leq&C_m(1-|a|^2)^{\delta/p}.\nonumber
\endr
 Similarly, when $j=1,2,\cdots,n$,
\begr
|K_m\circ (C_\varphi-C_\psi)(f_{a,j})(z)|\leq C_m(1-|a|^2)^{\delta/p},\nonumber
\endr
for some finite constant $C_m$ independent of $z$. Therefore, letting $i\to\infty$ and then using Fatou's Lemma as $m\to\infty$, by (\ref{1}), we have
$$\|C_\varphi-C_\psi-K\|_{A^p_\a\to A^q_\b}\gtrsim
\limsup_{i\to\infty}\|(C_\varphi-C_\psi)(f_{a_i,j})\|_{A^q_\b}.
$$
This remains valid for $1\leq q<\infty$ by a similar argument.
Therefore
\begr
\|C_\varphi-C_\psi\|_{e,A^p_\a\to A^q_\b}^q\geq C\limsup_{|a|\to1}\sum_{j=0}^n\|(C_\varphi-C_\psi)f_{a,j}\|_{A^q_\b}^q.\nonumber
\endr
The proof is complete.   $\Box$\msk

\noindent{\bf Proof of Theorem 1.2.}  First, we give the upper estimate for the essential norm of the   operator $C_\varphi-C_\psi$.
By the boundedness of the operator $C_\varphi-C_\psi:A^p_\a\to A^q_\b$, we have that $\omega_{\b,q}$ is a $(\lambda,\a)$-Bergman Carleson measure.
Defined by $K_n=C_{\varphi_n}$ for any $n\geq1$, where $\varphi_n(z)=\f{n}{n+1}z$. Every $K_n$ trivially has a norm less than $1$ and is compact on every $A^p_\a$ (\cite{Zhu2}). Let $R_n=I-K_n$. Then,
\begr
\|C_\varphi-C_\psi\|_{e,A^p_\a\to A^q_\a}&\leq&\limsup_{k\to\infty}\|(C_\varphi-C_\psi)R_k\|_{A^p_\a\to A^q_\b}\nonumber\\
&\leq&\limsup_{k\to\infty}\sup_{\|f\|_{A^p_\a}\leq1}\|(C_\varphi-C_\psi)R_kf\|_{A^q_\b}\nonumber.
\endr
Fix $s_0\in(0,1)$, set $E=\{z\in\B:\rho(z)\geq s_0\}$ and $E^\prime=\B\backslash E$. Let $f\in A^p_\a$ with $\|f\|_{A^p_\a}\leq1$.
Then we have 
\begr
I_k(f)&:=&\int_{E}|(C_\varphi-C_\psi)\circ R_kf(z)|^qd\nu_\b(z)\nonumber\\
&\leq&\left(\f{2}{s_0}\right)^q\int_{E}|R_kf(z)|^qd\omega_{\b,q}(z)\nonumber\\
&\lesssim&\int_{\B_s}|R_kf(z)|^qd\omega_{\b,q}(z)+\int_{\B\backslash \B_s}|R_kf(z)|^qd\omega_{\b,q}(z)\nonumber\\
&:=&I_{1,k}(f)+I_{2,k}(f).\nonumber
\endr
By Cauchy integral formula, we get
\begr
|R_k(f)(z)|&\leq& \f{1}{k+1}\sup_{w\in \B_s}|\Re f(w)|\lesssim\f{1}{k+1}\cdot\f{2}{1-s}\sup_{w\in\B_{\f{1+s}{2}}}|f(w)|\nonumber\\
&\leq&\f{2}{k(1-s)^{1+(n+1+\a)/p}}\lesssim \f{1}{k}\nonumber
\endr
for any $z\in \B_s$  with fixed $s\in (0,1)$. Here $\Re f$ denote the radial derivative of $f$. Therefore, for a fixed $s$, $\limsup\limits_{k\to\infty}I_{1,k}(f)=0.$

Now let us turn to $I_{2,k}$. We denote by $\omega_{\b,q}|_{\B\backslash\B_s}$ the restriction of $\omega_{\b,q}$ to
$\B\backslash \B_s$. Since $C_\varphi-C_\psi:A^p_\a\to A^q_\b$ is bounded by assumption, $\omega_{\b,q}$ is a $(\lambda,\a)$-Bergman Carleson measure, and $\omega_{\b,q}|_{\B\backslash\B_s}$ is also a $(\lambda,\a)$-Bergman Carleson measure.
Therefore,
\begr
I_{2,k}(f)&=&\|R_k(f)\|_{L^q(\B,\omega_{\b,q}|_{\B\backslash\B_s})}^q\lesssim\|\omega_{\b,q}|_{\B\backslash\B_s}\|_{\lambda,\a}\|R_k(f)\|_{A^p_\a}^q\nonumber\\
&\lesssim&\|\omega_{\b,q}|_{\B\backslash\B_s}\|_{\lambda,\a,r}.\nonumber
\endr
Letting $k\to\infty$ and $s\to 1$ in order, we obtain
$$\limsup_{k\to\infty}I_k(f)\lesssim\limsup_{s\to1}\|\omega_{\b,q}|_{\B\backslash\B_s}\|_{\lambda,\a,r}.$$

Denote
\begr
J_k(f)&=&\int_{E^\prime}|(C_\varphi-C_\psi)\circ R_kf(z)|^qd\nu_\b(z)\nonumber\\
&=&\left(\int_{E^\prime\cap\varphi^{-1}(\B_s)}+\int_{E^\prime\cap\varphi^{-1}(\B\backslash\B_s)}\right)|(C_\varphi-C_\psi)\circ R_kf(z)|^qd\nu_\b(z)\nonumber\\
&:=&J_{1,k}(f)+J_{2,k}(f).\nonumber
\endr
Let $r_0\in (0,1)$ be arbitrary. It is easy to see that
$$\lim_{k\to\infty}\sup_{\|f\|_{A^p_\a}\leq1}\sup_{|z|\leq r_0}|R_k(f)(z)|=0.$$
Thus,
$$\lim_{k\to\infty}\sup_{\|f\|_{A^p_\a}\leq1}\int_{E^\prime\cap\varphi^{-1}(\B_s)}|C_\varphi\circ R_k(f)(z)|^qd\nu_\b(z)=0$$
and
$$\lim_{k\to\infty}\sup_{\|f\|_{A^p_\a}\leq1}\int_{E^\prime\cap\varphi^{-1}(\B_s)}|C_\psi\circ R_k(f)(z)|^qd\nu_\b(z)=0,$$
here we used the fact that $E^\prime\cap\varphi^{-1}(\B_s)\subset\psi^{-1}(\B_{s^\prime})$, where $s^\prime=\f{s_0+s}{1+s_0s}$.
Hence
\begr
\limsup_{k\to\infty}J_k(f)
&\lesssim&\sup_{\|f\|_{A^p_\a}\leq1}\int_{F}|(C_\varphi-C_\psi)f(z)|^qd\nu_\b(z),\nonumber
\endr
where $F=E^\prime\cap\varphi^{-1}(\B\backslash\B_s)$.
In the estimate above we also used the fact that the operators $R_k$ are uniformly bounded.
Using Lemmas 2.1 and 2.2, Fubini's theorem and $1-|\varphi(z)|^2\asymp1-|w|^2$ for $\varphi(z)\in\triangle(w,r)$ we have
\begr
&&\int_{F}|(C_\varphi-C_\psi)f(z)|^qd\nu_\b(z)\nonumber\\
&\lesssim&\int_{F}\rho(z)^q\f{\int_{\triangle(\varphi(z),r)}|f(w)|^pd\nu_\a(w)}{(1-|\varphi(z)|^2)^{(n+1+\a)q/p}}d\nu_\b(z)\nonumber\\
&\lesssim&\int_{\B}|f(w)|^p\f{\int_{\varphi^{-1}(\triangle(w,r))\cap F}\rho(z)^qd\nu_\b(z)}{(1-|w|^2)^{(n+1+\a)q/p}}d\nu_\a(w)\nonumber\\
&\leq&\int_{\B}|f(w)|^p\f{\int_{\varphi^{-1}(\triangle(w,r)\cap (\B\backslash\B_s))}\rho(z)^qd\nu_\b(z)}{(1-|w|^2)^{(n+1+\a)q/p}}d\nu_\a(w)\nonumber\\
&\leq&\|f\|_{A^p_\a}^p\|\omega_{\b,q}|_{\B\backslash\B_s}\|_{\lambda,\a,r}.\nonumber
\endr
Letting $k\to\infty$ and $s\to1$ in order, and using the above estimate,  we have
\begr
\|C_\varphi-C_\psi\|_{e,A^p_\a\to A^q_\b}^q&\leq&\limsup_{k\to\infty}\|(C_\varphi-C_\psi)\circ R_k\|_{A^p_\a\to A^q_\b}^q\nonumber\\
&\lesssim&\limsup_{s\to1}\|\omega_{\b,q}|_{\B\backslash\B_s}\|_{\lambda,\a,r}.\nonumber
\endr

Next, we give the lower estimate for the essential norm of the   operator $C_\varphi-C_\psi$.
By Lemmas 2.4 and 3.1, we have
$$\|C_\varphi-C_\psi\|_{e,A^p_\a\to A^q_\b}\gtrsim \limsup_{|a|\to1}\int_{\varphi^{-1}(\triangle(a,r_0))}\f{\rho(z)^q}{(1-|a|^2)^{(n+1+\a)q/p}}d\nu_\b(z),$$
and
$$\|C_\varphi-C_\psi\|_{e,A^p_\a\to A^q_\b}\gtrsim \limsup_{|a|\to1}\int_{\psi^{-1}(\triangle(a,r_0))}\f{\rho(z)^q}{(1-|a|^2)^{(n+1+\a)q/p}}d\nu_\b(z).$$
Thus,
$$\|C_\varphi-C_\psi\|_{e,A^p_\a\to A^q_\b}\gtrsim \limsup_{|a|\to1}\f{\omega_{\b,q}(\triangle(a,r_0))}{(1-|a|^2)^{(n+1+\a)q/p}}.$$
The claim now follows by Lemma 2.5, (\ref{gs}) and (\ref{1.2}). The proof is complete.   $\Box$\msk

Let us now turn to the proof of the case $0<q<p<\infty$. For this we will make use of Khinchine's inequality. Define the Rademacher function $r_m$ by
$$r_m(t)=\textrm{sgn}(\sin(2^m\pi t)).$$
The Khinchine's inequality is the following.\msk

\noindent{\bf  Khinchine's inequality.}  {\it For $0<p<\infty$, there exist constants $0<A_p\leq B_p<\infty$ such that, for all natural numbers $m$ and all complex numbers $c_1, c_2, \cdots , c_m,$ we have
$$A_p\left(\sum_{j=1}^m |c_j|^2\right)^{\f{p}{2}}\leq\int_0^1\left|\sum_{i=1}^m c_jr_j(t)\right|^pdt\leq B_p\left(\sum_{j=1}^m |c_j|^2\right)^{\f{p}{2}}.$$
}\msk

\noindent{\bf Proof of Theorem 1.3.} We first prove that
$$\|C_\varphi-C_\psi\|_{A^p_\a\to A^q_\b}^q\lesssim\|\omega_{\b,q}\|_{\lambda,\a}.$$
Let $f\in A^p_\a$ with $\|f\|_{A^p_\a}\leq1$ and $r\in(0,1)$ be fixed. We write
 \begr
 && \|f\circ\varphi-f\circ\psi\|_{A^q_\b}^q\nonumber\\
 &=&  \Big(\int_{\{z\in\B:\rho(z)\geq r\}}+\int_{\{z\in\B:\rho(z)<r\}}\Big) |f\circ\varphi(z)-f\circ\psi(z)|^qd\nu_\b(z).\nonumber
 \endr
By the assumption on the joint pull-back measure $\omega_{\b,q}$, the first term is bounded above by
$\|\omega_{\b,q}\|_{\lambda,\a}$.
Applying Lemma 2.2, Fubini's Theorem and $1-|z|^2 \asymp 1-|w|^2$ for all $z\in\Delta(w,r)$ respectively, we see that the second term is bounded above by
\begr
&&C\int_{\B}|f(w)|^q\f{\int_{\varphi^{-1}(\triangle(w,r^\prime))\cap \{z\in\B:\rho(z)<r\}}\rho(z)^qd\nu_\b(z)}{(1-|w|^2)^{n+1+\a}}d\nu_\a(w)\nonumber\\
&\leq& C\int_{\B}|f(w)|^q\f{\omega_{\b,q}(\triangle(w,r^\prime))}{(1-|w|^2)^{n+1+\a}}d\nu_\a(w),\nonumber
\endr
where the constants $C$ and $r^\prime\in(0,1)$ depend only on $n$, $\a$ and $r$. Applying the H\"{o}lder's inequality and Theorem C we get the desired result.

Next, we prove that
$$\|\{\omega_{\b,q,k}\}\|_{l^{\f{p}{p-q}}}\lesssim\|C_\varphi-C_\psi\|_{A^p_\a\to A^q_\b}.$$
Here
$$\omega_{\b,q,k}=\f{\omega_{\b,q}(\triangle(a_k,r_0))}{(1-|a_k|)^{(n+1+\a)q/p}}.
$$ We suppose $M:=\|C_\varphi-C_\psi\|_{A^p_\a\to A^q_\b}<\infty$. If $M=0$, it is easy to see that $\varphi=\psi$ and $\omega_{\b,q}=0$ and the claim is straightforward. So, we suppose $M\neq0$.
Let $0<r_0<1$ and $N=N(r_0)$ be defined as in Lemma 2.3 and denote $r_1=1-\f{1}{2N}$.
Let $\{a_k\}$ be an $r_0$-lattice of $\B$ in the pseudo-hyperbolic metric
with $|a_1|\leq |a_2|\leq\cdots\leq |a_n|\leq\cdots$.
By Lemma 2.6, there exist a nonnegative number $$L_0\leq\f{(\f{2}{r_0}+1)^{2n}}{(1-r_1^2)^n}$$ such that
$|a_1|\leq |a_2|\leq\cdots\leq|a_{L_0}|\leq r_1$ and $|a_k|>r_1$ for any $k\geq L_0+1$.
For $\{c_j\}\in l^p$,
define
$$g_j(z)=\sum\limits_{k=1}^\infty c_kf_{a_k,j}(z),j=1,2,\cdot\cdot\cdot, n,$$
where $f_{a,j}$ are defined as Lemma 2.4.
Then $\|g_{j}\|_{A^p_\a}\lesssim(\sum\limits_{k=1}^\infty|c_k|^p)^{\f{1}{p}}.$
Using the boundedness of the operator $C_\varphi-C_\psi$, we have
\begr
M^q\Big (\sum_{k=1}^\infty|c_k|^p\Big )^{\f{q}{p}}&\geq& C\|(C_\varphi-C_\psi)g_j\|_{A^q_\b}^q\nonumber\\
&=&C\int_{\B}|\sum_{k=1}^\infty c_k(C_\varphi-C_\psi)\circ f_{a_k,j}(z)|^qd\nu_\b(z).\nonumber
\endr
Replace now $c_k$ by $r_k(t)c_k$ and then integrate with respect to $t$ from $0$ to $1$. Applying Fubini's theorem and the Khinchine's inequality, we get
\begr
M^q\Big (\sum_{k=1}^\infty|c_k|^p\Big )^{\f{q}{p}}&\gtrsim& \int_{\B}\int_0^1|\sum_{k=1}^\infty c_kr_k(t)(C_\varphi-C_\psi)\circ f_{a_k,j}(z)|^qdtd\nu_\b(z)\nonumber\\
&\gtrsim& A_q\int_{\B}\left(\sum_{k=1}^\infty |c_k|^2|f_{a_k,j}\circ\varphi(z)-f_{a_k,j}\circ\psi(z)|^2\right)^{\f{q}{2}}d\nu_\b(z).\nonumber
\endr
Applying Lemma 2.4 and $|1-\langle z,a_k\rangle|\asymp 1-|a_k|^2$ for $z\in\triangle(a_k,r_0)$, we have
\begr
&&\sum_{j=0}^n |f_{a_k,j}\circ\varphi(z)-f_{a_k,j}\circ\psi(z)|^2\nonumber\\
&\gtrsim&|f_{a_k,0}\circ\varphi(z)|^2\rho(z)^2\chi_{\varphi^{-1}(\triangle(a_k,r_0))}(z)\nonumber\\
&\gtrsim&\f{\rho(z)^2\chi_{\varphi^{-1}(\triangle(a_k,r_0))}(z)}{(1-|a_k|)^{(n+1+\a)2/p}}.\nonumber
\endr
By Lemma 2.6, there exists a $$K=K(r_0)\leq\f{(\f{2}{r_0}+1)^{2n}}{(1-r_0^2)^n}$$
  such that
for any $z\in \B$, there are at most $K$ elements of $\{a_k\}$ lying in $\triangle(z,r_0)$.  Therefore
\begr
&&(n+1)M^q\Big (\sum_{k=1}^\infty|c_k|^p\Big )^{\f{q}{p}}\nonumber\\
&\gtrsim& A_q\max\{1,(n+1)^{1-\f{q}{2}}\}\nonumber\\
&&\cdot\int_{\B}\left(\sum_{k=1}^\infty |c_k|^2\sum_{j=0}^n|f_{a_k,j}\circ\varphi(z)-f_{a_k,j}\circ\psi(z)|^2\right)^{\f{q}{2}}d\nu_\b(z)\nonumber\\
&\gtrsim&\int_{\B}\left(\sum_{k=L_0+1}^\infty |c_k|^2\f{\rho(z)^2\chi_{\varphi^{-1}(\triangle(a_k,r_0))}(z)}{(1-|a_k|)^{(n+1+\a)2/p}}\right)^{\f{q}{2}}d\nu_\b(z)\nonumber\\
&\geq&\max\{1,K^{\f{q}{2}-1}\}\int_{\B}\sum_{k=L_0+1}^\infty |c_k|^q\f{\rho(z)^q\chi_{\varphi^{-1}(\triangle(a_k,r_0))}(z)}{(1-|a_k|)^{(n+1+\a)q/p}}d\nu_\b(z)\nonumber\\
&\gtrsim&\sum_{k=L_0+1}^\infty |c_k|^q\f{\int_{\varphi^{-1}(\triangle(a_k,r_0))}\rho(z)^qd\nu_\b(z)}{(1-|a_k|)^{(n+1+\a)q/p}}.\nonumber
\endr
Change the roles of $\varphi$ and $\psi$, we get
\begr
(n+1)M^q\Big(\sum_{k=1}^\infty|c_k|^p\Big)^{\f{q}{p}}
&\gtrsim&\sum_{k=L_0+1}^\infty |c_k|^q\f{\int_{\psi^{-1}(\triangle(a_k,r_0))}\rho(z)^qd\nu_\b(z)}{(1-|a_k|)^{(n+1+\a)q/p}}.\nonumber
\endr
Therefore
\begr
M^q\Big(\sum_{k=1}^\infty|c_k|^p\Big)^{\f{q}{p}}
&\gtrsim&\sum_{k=L_0+1}^\infty |c_k|^q\f{\omega_{\b,q}(\triangle(a_k,r_0))}{(1-|a_k|)^{(n+1+\a)q/p}}.\nonumber
\endr
Let
$$b_k=\f{1}{M^q}\f{\omega_{\b,q}(\triangle(a_k,r_0))}{(1-|a_k|^2)^{(n+1+\a)q/p}}, k=1,2,3,\cdots.$$
Since $\{c_k\}\in l^p$ is arbitrary, we deduce
$$(b_{L_0+1},b_{L_0+2},\cdots)\in (l^{\f{p}{q}})^{\star}=l^{\f{p}{p-q}}$$
and there exists a constant $C>0$ such that
\begr
\sum_{k=L_0+1}^\infty \left(\f{\omega_{\b,q}(\triangle(a_k,r_0))}{(1-|a_k|)^{(n+1+\a)q/p}}\right)^{\f{p}{p-q}}<CM^{\f{pq}{p-q}} .\nonumber
\endr
For $1\leq i\leq L_0$, we have $|a_i|\leq 1-\f{1}{2N}$.  Denote $r_1=1-\f{1}{2N}$ and $r_2=\f{r_0+r_1}{1+r_0r_1}$. Then $\triangle(a_i,r_0)\subset\triangle(0,r_2)$. By Lemma 2.1, we have
\begr
b_i&\lesssim& \f{1}{M^q(2N)^{\frac{(n+1+\a)q}{p}}}\left(\int_{\varphi^{-1}(\triangle(a_i,r_0))} +\int_{\psi^{-1}(\triangle(a_i,r_0))}\right)\rho(z)^qd\nu_\b(z)\nonumber\\
&\leq& \f{1}{M^q(2N)^{\frac{(n+1+\a)q}{p}}}\left(\int_{\varphi^{-1}(\triangle(a_i,r_0))}+\int_{\psi^{-1}(\triangle(a_i,r_0))}\right)\f{|\varphi(z)-\psi(z)|^q}{(1-r_2)^q}d\nu_\b(z)\nonumber\\
&\leq& \f{1}{M^q(2N)^{\frac{(n+1+\a)q}{p}}}\cdot\f{2M^q\|z\|_{A^p_\a}^q}{(1-r_2)^q}\nonumber\\
&\leq&C(r,N,n,p,q,\a).  \nonumber
\endr
 Thus, we get
$$\{ b_1,b_{2},\cdots\}\in l^{\f{p}{p-q}}$$
and there exist a constant $C>0$ such that
\begr
\sum_{k=1}^\infty \left(\f{\omega_{\b,q}(\triangle(a_k,r_0))}{(1-|a_k|)^{(n+1+\a)q/p}}\right)^{\f{p}{p-q}}<CM^{\f{pq}{p-q}} .\nonumber
\endr
Hence,
 $\omega_{\b,q}$ is a $(\lambda,\a)$-Bergman Carleson measure, and
$$\|\{\omega_{\b,q,k}\}\|_{l^{p/(p-q)}}\lesssim\|C_\varphi-C_\psi\|_{A^p_\a\to A^q_\b}^q.$$
Therefore, by the above discussion and Theorem C, we get
\begr
\|C_\varphi-C_\psi\|_{A^p_\a\to A^q_\b}&\asymp&\|\omega_{\b,q}\|_{\lambda,\a}\asymp \|\{\omega_{\b,q,k}\} \|_{l^{p/(p-q)}}\nonumber\\
&\asymp&\|\widehat{\omega_{\b,q}}_r\|_{L^{p/(p-q)}_{\nu_\a}}
\asymp\|B_{s,\a}(\omega_{\b,q})\|_{L^{p/(p-q)}_{\nu_\a}}.\nonumber
\endr

Finally, we suppose that $\omega_{\b,q}$ is a $(\lambda,\a)$-Bergman Carleson measure, and then prove that the operator $C_\varphi-C_\psi:A^p_\a\to A^q_\b$ is compact.
Let $\{f_k\}$ be any sequence in
$A^p_\a$ with $\|f_k\|_{A^p_\a}\leq1$ and ${f_k}\to0$ uniformly on  compact subsets of $\B$.
Then by the Remark after Theorem C, $\omega_{\b,q}$ is a vanishing $(\lambda,\a)$-Bergman Carleson measure. Therefore, we have
$$\lim_{k\to\infty}\int_{\B}|f_k(z)|^qd\omega_{\b,q}(z)=0.$$
For $r\in(0,1)$ fixed, we write
\begr &&\|f_k\circ\varphi-f_k\circ\psi\|_{A^q_\b}^q\nonumber\\
&=&\Big(\int_{\{z\in\B:\rho(z)\geq r\}}+\int_{\{z\in\B:\rho(z)<r\}}\Big)|f_k\circ\varphi(z)-f_k\circ\psi(z)|^qd\nu_\b(z).\nonumber
\endr
The fist term is uniformly bounded above by
$$\Big(\f{2}{r}\Big)^q\int_{\{z\in\B:\rho(z)\geq r\}}|f_k(z)|^qd\omega_{\b,q}(z).$$%\to0,\mbox{~~~as~~}k\to\infty.$$
For any fixed $s\in(0,1)$,  we see that the second term is bounded above by
\begr
&&C\int_{\B}|f_k(w)|^q\f{\omega_{\b,q}(\triangle(w,r^\prime))}{(1-|w|^2)^{n+1+\a}}d\nu_\a(w)\nonumber\\
&=&C\left(\int_{\B_s}+\int_{\B\backslash\B_s}\right)|f_k(w)|^q\f{\omega_{\b,q}(\triangle(w,r^\prime))}{(1-|w|^2)^{n+1+\a}}d\nu_\a(w)\nonumber\\
&:=&I_1+I_2,\nonumber
\endr
where the constants $C$ and $r^\prime\in (0,1)$ depend only on $n$, $\a$ and $r$.
Since $\omega_{\b,q}(\B)\leq2$, we have
$$I_1\lesssim\int_{\B_s}|f_k(z)|^qd\nu_\a(z).$$
Applying the H\"{o}lder's inequality, we have
\begr
I_2&\lesssim&\|f_k\|_{A^p_\a}^q\left\|\f{\omega_{\b,q}(\triangle(w,r^\prime))}{(1-|w|^2)^{n+1+\a}}\right\|_{L^t(\B\backslash\B_s,\omega_{\b,q})}\nonumber\\
&\leq&\left\|\f{\omega_{\b,q}(\triangle(w,r^\prime))}{(1-|w|^2)^{n+1+\a}}\right\|_{L^t(\B\backslash\B_s,\omega_{\b,q})}.\nonumber
\endr
Letting $k\to\infty$, we obtain
$$\limsup_{k\to\infty}\|f_k\circ\varphi-f_k\circ\psi\|_{A^q_\b}^q
\lesssim\left\|\f{\omega_{\b,q}(\triangle(w,r^\prime))}{(1-|w|^2)^{n+1+\a}}\right\|_{L^t(\B\backslash\B_s,\omega_{\b,q})}.$$
Since $\omega_{\b,q}$ is a $(\lambda,\a)$-Bergman Carleson measure, we have
$$\f{\omega_{\b,q}(\triangle(w,r^\prime))}{(1-|w|^2)^{n+1+\a}}\in L^t(\B,\omega_{\b,q}).$$
Thus, letting $s\to1$, by the Lebesgue Dominated Covergence Theorem, we get
$$\lim_{k\to\infty}\|f_k\circ\varphi-f_k\circ\psi\|_{A^q_\b}=0.$$
Therefore, $C_\varphi-C_\psi:A^p_\a\to A^q_\b$ is compact. The proof is complete.         $\square$ \msk

\noindent{\bf Remark.}  {The methods we used to prove Theorem 1.1-1.3 can be generalized to the case when $\nu_\b$ is replaced by a positive Borel measure $\mu$ by using the same technique.}

\end{document}